\documentclass[12pt,oneside]{article}
\usepackage{amsmath,amssymb,amsfonts,amsthm}
\usepackage{color}
\usepackage[all]{xy}
\usepackage{graphicx}
\usepackage{color}
\usepackage{makeidx}
\usepackage{multirow}
\usepackage{hyperref}
\usepackage{rotating}

\allowdisplaybreaks

\def\paa{\dot{\partial}}

\numberwithin{equation}{section} 
\numberwithin{figure}{section} 

\theoremstyle{plain}
\newtheorem*{theorem*}{Theorem}
\newtheorem{theorem}{Theorem}[section]
\newtheorem{lemma}[theorem]{Lemma}

\newtheorem{corollary}[theorem]{Corollary}

\theoremstyle{definition}
\newtheorem{definition}[theorem]{Definition}

\theoremstyle{remark}

\newtheorem*{acknowledgement*}{Acknowledgement}

\numberwithin{equation}{section}

\newcommand\overcirc[1]{\raisebox{10pt}{\tiny{$\circ$}}{\kern-7.5pt}\mbox{$#1$}}
\newcommand\undersym[2]{\raisebox{-6pt}{$#2$}{\kern-5pt}\mbox{$#1$}}
\newcommand\overdiamond[1]{\raisebox{10pt}{\small$\star$}{\kern-7.5pt}\mbox{$#1$}}
\newcommand\overast[1]{\raisebox{10pt}{\small$\ast$}{\kern-7.5pt}\mbox{$#1$}}
\newcommand\overlind[1]{\raisebox{10pt}{\small$\overline{{\hspace{2pt}}\star}$}{\kern-7.5pt}\mbox{$#1$}}
\newcommand\overlinc[1]{\raisebox{10pt}{\small$\overline{{\hspace{2pt}}\circ}$}{\kern-7.5pt}\mbox{$#1$}}
\newcommand\overlina[1]{\raisebox{10pt}{\small$\overline{{\hspace{1pt}}\ast}$}{\kern-7.5pt}\mbox{$#1$}}

\begin{document}

\title{\bf  On Landsberg Unicorns of General $(\alpha,\beta)$-Type  
 }

 \author{Salah G. Elgendi}

\date{}

\maketitle
\vspace{-1.20cm}

\begin{center}
   Department of Mathematics, Faculty of Science,
  \\
  Islamic University of Madinah, Madinah, Saudi Arabia
\end{center}

\begin{center}
  salah.ali@fsc.bu.edu.eg, salahelgendi@yahoo.com
\end{center}

\vspace{0.3cm}

  \begin{abstract}
In this paper, we investigate the classes of \lq\lq unicorns" of general $(\alpha,\beta)$-metric types, namely, Landsberg general $(\alpha,\beta)$-metrics that are non-Berwaldian. We provide explicit solutions for these metrics, including several new classes. Furthermore, we identify and   exclude the cases that yield degenerate metrics to guarantee the geometric validity of the established solutions.
\end{abstract}



\section{Introduction} 
~\par

In Finsler geometry, one of the most celebrated and long-standing open problems is the Landsberg-Berwald puzzle: \textit{Does there exist a regular Landsberg metric that is not a Berwald metric?} 

To recall the underlying geometric distinction, a Finsler manifold $(M,F)$ is called a  \textit{{Berwald space}} if its Berwald connection coefficients depend solely on the position coordinates, which implies that its geodesic spray is quadratic. Conversely, it is called a \textit{Landsberg space} if the metric tensor is horizontally covariantly constant with respect to the Berwald connection. While every Berwald metric is inherently Landsbergian, finding a regular converse example, a metric satisfying the latter property but not the former, has remained so elusive that David Bao famously dubbed these hypothetical spaces the \textit{``unicorns'' of Finsler geometry}~\cite{Bao07}. More information and surveys on the problem can be found in \cite{Crampin_2,Tayebi2021}.

To search for these mathematical entities, geometers frequently turn to the class of {$(\alpha,\beta)$-metrics}. Defined in the form $F = \alpha \phi(\beta/\alpha)$, where $\alpha = \sqrt{a_{ij}(x)y^iy^j}$ is a Riemannian metric and $\beta = b_i(x)y^i$ is a 1-form, $(\alpha,\beta)$-metrics provide a highly computable yet geometrically rich testing ground. In a seminal paper, Zhongmin Shen demonstrated that there are no regular, non-Berwaldian Landsberg $(\alpha,\beta)$-metrics~\cite{Shen2009}. However, if the regularity condition is relaxed by allowing for structural singularities,  the mathematical landscape changes completely. Asanov successfully broke ground by uncovering non-regular ``unicorns'' within specific subclasses~\cite{Asa06}, sparking an active hunt for explicit, non-regular Landsberg metrics that fail to reduce to Berwald types.

In this direction, Elgendi provided several explicit examples of such unicorns within the frameworks of both $(\alpha,\beta)$-metrics and spherically symmetric metrics (see \cite{Elgendi-solutions, Elgendi-JGP, Elgendi-CM, Elgendi-MJM} and references therein). Further significant progress was achieved by Zhou, Li and Wang  \cite{ZhouLi2018, Zhou-Wang-Li}, who classified the unicorns of  general $(\alpha,\beta)$-types.

\medskip

Motivated by these developments, this paper focuses on the classes of non-regular non-Berwaldian Landsberg metrics, of general $(\alpha,\beta)$-types, introduced in \cite{Zhou-Wang-Li}. For the partial differential equations (PDEs) that characterize the metric function $\phi$, we explicitly compute the corresponding integrals and provide complete analytical formulas for $\phi(s, b^2)$. Furthermore, we unveil new families of examples and establish rigorous parametric restrictions on the arbitrary functions and constants involved, ensuring that the resulting metrics strictly avoid degeneracy.

The primary contribution of this work is the complete classification and explicit determination of the metric functions $\phi$ for  non-regular (including almost regular) non-Berwaldian Landsberg general $(\alpha,\beta)$-metrics. Our analysis yields three distinct structural classes of solutions. Crucially, each class naturally bifurcates into three distinct cases depending on the sign of an algebraic quantity: \textit{hyperbolic}, \textit{trigonometric}, and \textit{rational}.

\medskip

In what follows, we present the explicit forms of these non-regular unicorns of general $(\alpha,\beta)$-metric type and of dimension $n\geq3$. Throughout this paper, we use the term non-regular rather than almost regular,  because the former constitutes a broader class of metrics.  For most of the following classes, non-regularity arises precisely when $s=\pm b$. 

\subsection*{Class I Metrics}
The first family of metrics arises from the structural conditions (see Theorem~\ref{Thm:Class_I}). Governed by the quantity $\Delta_1 = c_1^2 - 4c_2 - 4$, the explicit analytical form of the metric function $\phi(s)$ is completely determined by the following piecewise representation:
\begin{equation*}
\label{eq:phi_class_I_main} \phi(s) = \begin{cases} c_3 \sqrt{b^2 + c_2 s^2 + c_1 s \sqrt{b^2-s^2}} \exp\left( \frac{c_1 }{\sqrt{\Delta_1}} \operatorname{arctanh} \left( \frac{c_1 s + 2\sqrt{b^2-s^2}}{s \sqrt{\Delta_1}} \right) \right), & \text{if } \Delta_1 > 0, \\ c_3 \sqrt{b^2 + c_2 s^2 + c_1 s \sqrt{b^2-s^2}} \exp\left(- \frac{c_1 }{\sqrt{-\Delta_1}} \operatorname{arctan} \left( \frac{c_1 s + 2\sqrt{b^2-s^2}}{s \sqrt{-\Delta_1}} \right) \right), & \text{if } \Delta_1 < 0, \\ c_3 \left( c_1 s + 2\sqrt{b^2-s^2} \right) \exp\left( \frac{c_1 s}{c_1 s + 2\sqrt{b^2-s^2}} \right), & \text{if } \Delta_1 = 0, \end{cases}
\end{equation*}
where $c_1 \neq 0$ is the non-Berwaldian condition, $c_2 \neq -1$ ensures non-degeneracy, and $c_3 \neq 0$ is a real constant. Furthermore, the 1-form $\beta$ satisfies
$$b = b_0, \quad b_{i \mid j} = k(x) \left(b^2 a_{i j} - b_i b_j\right), \quad k(x)\neq 0.$$ 
The third class first appeared in a special case in \cite{Elgendi-solutions}.   Subsequently, Heefer et al. \cite{Heefer2023}  employed it to construct a new family of exact vacuum solutions  to Pfeifer and Wohlfarth’s field equations in Finsler gravity.

\subsection*{Class II Metrics}
The second family represents general $(\alpha,\beta)$-metrics  (see Theorem~\ref{Thm:Class_II}). Here, we have the quantity $\Delta_2 = c_1^2 -(\tilde{b}^2\lambda_1 + 1)(1 - 2\tilde{b}^2\lambda_0)$, yielding the  result:
\begin{equation*}\label{eq:phi_class_II_main} \phi(\tilde{b}^2, \tilde{s}) = \begin{cases} \lambda_3 \sqrt{2 \tilde{s} c_1 \sqrt{\tilde{b}^2-\tilde{s}^2} + \tilde{b}^2 \tilde{s}^2 \lambda_1 + \tilde{b}^2 - 2\tilde{b}^2(\tilde{b}^2-\tilde{s}^2)\lambda_0}\, e^{ -\frac{c_1}{\sqrt{\Delta_2}} \operatorname{arctanh} \left( \frac{(\tilde{b}^2\lambda_1+1)\tilde{s} + c_1\sqrt{\tilde{b}^2-\tilde{s}^2}}{\sqrt{\Delta_2}\sqrt{\tilde{b}^2-\tilde{s}^2}} \right) }, & \text{if } \Delta_2 > 0, \\ \lambda_3 \sqrt{2 \tilde{s} c_1 \sqrt{\tilde{b}^2-\tilde{s}^2} + \tilde{b}^2 \tilde{s}^2 \lambda_1 + \tilde{b}^2 - 2\tilde{b}^2(\tilde{b}^2-\tilde{s}^2)\lambda_0}\, e^{ \frac{c_1}{\sqrt{-\Delta_2}} \operatorname{arctan} \left( \frac{(\tilde{b}^2\lambda_1+1)\tilde{s} + c_1\sqrt{\tilde{b}^2-\tilde{s}^2}}{\sqrt{-\Delta_2}\sqrt{\tilde{b}^2-\tilde{s}^2}} \right) }, & \text{if } \Delta_2 < 0, \\ \lambda_3 \left( (\tilde{b}^2\lambda_1+1)\tilde{s} + c_1\sqrt{\tilde{b}^2-\tilde{s}^2} \right) e^{ - \frac{c_1 \sqrt{\tilde{b}^2-\tilde{s}^2}}{(\tilde{b}^2\lambda_1+1)\tilde{s} + c_1\sqrt{\tilde{b}^2-\tilde{s}^2}} }, & \text{if } \Delta_2 = 0, \end{cases}
\end{equation*}
where $\lambda_3(\tilde{b}^2)$ is a non-zero smooth function, $1-2\tilde{b}^2\lambda_0 \neq 0$, and $\tilde{b}^2\lambda_1+1 \neq 0$. The 1-form $\beta$ satisfies
$$\widetilde{b}_{i \mid j}=\widetilde{c} a_{i j}.$$ 

\subsection*{Class III/IV Metrics}
Finally, Class III/IV metrics correspond to solutions of  the PDE that determines the $\phi$ (Theorem~\ref{Thm:Class_III_IV}). Controlled by the quantity $\Delta_3 = b^4 f_1(b^2)^2 - 4f_2(b^2)$, the function $\phi$ is given by:
\begin{equation*}\label{eq:phi_class_III_main} \phi(b^2, s) = \begin{cases} f_3 \sqrt{b^2 + b^2 s \sqrt{b^2-s^2} f_1(b^2) + s^2(f_2(b^2) - 1)} \, e^{ -\frac{b^2 f_1(b^2)}{ \sqrt{\Delta_3}} \operatorname{arctanh} \left( \frac{b^2 f_1(b^2) \sqrt{b^2-s^2} + 2s f_2(b^2)}{\sqrt{\Delta_3} \sqrt{b^2-s^2}} \right) }, & \text{if } \Delta_3 > 0, \\ 
f_3 \sqrt{b^2 + b^2 s \sqrt{b^2-s^2} f_1(b^2) + s^2(f_2(b^2) - 1)} \, e^{ \frac{b^2 f_1(b^2)}{ \sqrt{-\Delta_3}} \operatorname{arctan} \left( \frac{b^2 f_1(b^2) \sqrt{b^2-s^2} + 2s f_2(b^2)}{\sqrt{-\Delta_3} \sqrt{b^2-s^2}} \right) }, & \text{if } \Delta_3 < 0, \\
 f_3 \left( b^2 f_1(b^2) s + 2 \sqrt{b^2-s^2} \right) e^{ -\frac{2 \sqrt{b^2-s^2} }{ b^2 f_1(b^2) s + 2\sqrt{b^2-s^2}} }, & \text{if } \Delta_3 = 0, \end{cases} 
\end{equation*}
where $f_3(b^2) > 0$ is an arbitrary smooth  function, and $f_2(b^2) = C b^4 f_1(b^2)^2$ for some real constant $C$. In this setup, the 1-form $\beta$ satisfies either
\begin{align*}
s_{i j} &= \frac{1}{b^2}\left(b_i s_j-b_j s_i\right), \\ 
r_{i j} &= k\left(b^2 a_{i j}-b_i b_j\right)+c b_i b_j+\lambda\left(b_i s_j+b_j s_i\right), 
\end{align*}
where $k=k(x)$, $c=c(x)$, and $\lambda=\lambda(x) \neq -\frac{1}{b^2}$ are scalar functions, or alternatively,
\begin{align*} 
b_{i \mid j}=k\left(b^2 a_{i j}-b_i b_j\right)+c b_i b_j+\theta_i b_j+\theta_j b_i, 
\end{align*}
where $\theta_i$ is a 1-form perpendicular to $\beta$.

By establishing these explicit piecewise representations at the outset, we obtain a direct correspondence between the underlying algebraic constraints and the resulting analytical forms of  function $\phi$.

\section{Notations and preliminaries}
~\par

Let $M$ be an $n$-dimensional smooth manifold. The coordinates on $M$ are denoted by $x^i$ and the corresponding fibre coordinates by $y^i$. The tangent bundle of $M$ is denoted by $TM$, while the slit tangent bundle is denoted by $T M \setminus \{0\}$.

We denote partial derivatives by $\partial_i := \frac{\partial}{\partial x^i}$ and $\dot{\partial}_i := \frac{\partial}{\partial y^i}$.

\begin{definition}
A \emph{Finsler manifold} of dimension $n$ is a pair $(M, F)$, where $M$ is an $n$-dimensional smooth manifold and $F : TM \rightarrow \mathbb{R}$ is a function satisfying:
\begin{itemize}
  \item[(a)] $F$ is smooth and positive on $T M \setminus \{0\}$ and vanishes only at the zero section,
  \item[(b)] $F$ is positively homogeneous of degree one in $y$,
  \item[(c)] The metric tensor components $g_{ij} = \frac{1}{2}\dot{\partial}_i \dot{\partial}_j F^2$ form a non-degenerate matrix $(g_{ij})$ on $T M \setminus \{0\}$.
\end{itemize}
Such a pair $(M, F)$ is called  a \emph{regular Finsler space}. If $F$ fails to be smooth or defined in certain directions, it is termed a \emph{non-regular} Finsler function. Under certain conditions, $(M,F)$ is called a conic Finsler manifold when $F$ fulfills the above conditions on a conic subset of $T M \setminus \{0\}$. 
\end{definition}

The geodesics of a Finsler metric $F=F(x, y)$ on $M$ are characterized by the system of second-order ordinary differential equations:
\begin{equation}
\frac{d^{2} x^{i}}{d t^{2}}+2 G^{i}\left(x, \frac{d x}{d t}\right)=0,
\end{equation}
where the geodesic (or spray) coefficients $G^{i}=G^{i}(x, y)$ are defined by
\begin{equation}
G^{i} := \frac{1}{4} g^{i l} \left\{ y^m \partial_m \dot{\partial}_l F^2 - \partial_l F^2 \right\}.
\end{equation}
Here, $(g^{i j})$ represents the inverse matrix of $(g_{ij})$. 

The \textit{Berwald tensor} $B = B_{j k \ell}^{i} \partial_{i} \otimes d x^{j} \otimes d x^{k} \otimes d x^{\ell}$ is generated by taking vertical derivatives of the spray coefficients:
\begin{equation}
B_{j k \ell}^{i} := \paa_j\paa_k\paa_\ell G^{i}. 
\end{equation}
A Finsler metric is classified as a \textit{Berwald metric} if $B_{j k \ell}^{i} = 0$, which geometrically means that the coefficients $G^{i}(x, y)$ depend quadratically on $y \in T_{x} M$ for every point $x \in M$.

Similarly, the \textit{Landsberg tensor} $L = L_{j k \ell} d x^{j} \otimes d x^{k} \otimes d x^{\ell}$ is defined via
\begin{equation}
L_{j k \ell} := -\frac{1}{2} F B_{j k \ell}^{i} \paa_i F  . 
\end{equation}
Metrics satisfying $L_{j k \ell} = 0$ are called \textit{Landsberg metrics}. By definition, every Berwald metric is inherently a Landsberg metric.

In 2011, a comprehensive metric family called a general $(\alpha, \beta)$-metric was introduced by Yu and Zhu \cite{Yu-Zhu} as follows:
$$
F=\alpha \phi\left(b^2, s\right), \quad s=\frac{\beta}{\alpha},
$$
where $\phi(s)$ is a $C^{\infty}$ function, $\alpha=\sqrt{a_{i j}(x) y^i y^j}$ is a Riemannian metric, $\beta=b_i(x) y^i$ is a 1-form, and $b:=\left\|\beta_x\right\|_\alpha$. This structural framework directly generalizes the classic class of $(\alpha, \beta)$-metrics.

To guarantee the positive definiteness of general $(\alpha, \beta)$-metrics, we recall the following regularity criterion:
\begin{lemma}\cite{Yu-Zhu}
Let $M$ be an $n$-dimensional manifold. The function $F = \alpha \phi\left(b^{2}, \frac{\beta}{\alpha}\right)$ constitutes a valid Finsler metric on $M$ for any Riemannian metric $\alpha$ and $1$-form $\beta$ with length $\|\beta\|_{\alpha} < b_{0}$ if and only if $\phi = \phi\left(b^{2}, s\right)$ is a positive smooth function satisfying 
\begin{equation}
\phi - s \phi_{s} > 0, \quad \phi - s \phi_{s} + \left(b^{2}-s^{2}\right) \phi_{ss} > 0 
\end{equation}
\textit{for $n \geqslant 3$, or simply $\phi - s \phi_{s} + \left(b^{2}-s^{2}\right) \phi_{ss} > 0$ when $n = 2$, across all variables bounded by $|s| \leqslant b < b_{0}$, where the subscript $s$ denotes the partial derivative with respect to $s$.}
\end{lemma}
The determinant $\operatorname{det}\left(g_{i j}\right)$ is given by \cite{Yu-Zhu}:
$$\operatorname{det}\left(g_{i j}\right)= \phi^{n+1}\left(\phi-s \phi_s\right)^{n-2}\left(\phi-s \phi_s+\left(b^2-s^2\right) \phi_{ss}\right) \operatorname{det}\left(a_{i j}\right),$$
where $\operatorname{det}\left(a_{i j}\right)$ is the determinant of the Riemannian metric $\alpha$. 

Depending on this formula of $\operatorname{det}\left(g_{i j}\right)$, we can provide the non-degeneracy conditions of the resulting Finsler metric. In this work, it is appropriate to write the non-degeneracy conditions as follows:
\begin{align}
\label{Eq:simple_form_almost regularity_1} \phi - s \phi_s &= \phi\left( 1 - s\frac{\phi_s}{\phi} \right) \neq 0, \\
\label{Eq:simple_form_almost regularity_2} \phi - s \phi_s + (b^2 - s^2)\phi_{ss} &= \phi\left( 1 - s\frac{\phi_s}{\phi} + (b^2-s^2) \left[ \left(\frac{\phi_s}{\phi}\right)_s + \left(\frac{\phi_s}{\phi}\right)^2 \right] \right) \neq 0.
\end{align}
If either factor vanishes, then $\det(g_{ij})=0$ and the metric tensor becomes degenerate. 

\medskip

Let $b_{i \mid j}$ denote the covariant derivative of $\beta$ with respect to the Levi-Civita connection of $\alpha$. We employ the following standard abbreviations:
$$
\begin{aligned}
& r_{i j} := \frac{1}{2}(b_{i \mid j}+b_{j \mid i}), \quad s_{i j} := \frac{1}{2}(b_{i \mid j}-b_{j \mid i}), \quad r_{00} := r_{i j} y^{i} y^{j}, \quad s_{0}^{i} := a^{i j} s_{j k} y^{k}, \\
& r_{i} := b^{j} r_{j i}, \quad s_{i} := b^{j} s_{j i}, \quad r_{0} := r_{i} y^{i}, \quad s_{0} := s_{i} y^{i}, \quad r^{i} := a^{i j} r_{j}, \quad s^{i} := a^{i j} s_{j}, \quad r := b^{i} r_{i}.
\end{aligned}
$$

According to \cite{Yu-Zhu}, the spray coefficients $G^{i}$ of a general $(\alpha, \beta)$-metric can be expressed as:
\begin{equation}
G^{i} = G_{\alpha}^{i} + P y^{i} + Q^{i}, 
\end{equation}
where $G_{\alpha}^{i}$ are the geodesic coefficients of $\alpha$, and the remaining terms are given by:
$$
\begin{aligned}
P &= \left\{\Theta\left(-2 \alpha Q s_{0}+r_{00}+2 \alpha^{2} R r\right)+\alpha \Omega\left(r_{0}+s_{0}\right)\right\} \alpha^{-1}, \\
Q^{i} &= \alpha Q s_{0}^{i}-\alpha^{2} R\left(r^{i}+s^{i}\right)+\left\{\Psi\left(-2 \alpha Q s_{0}+r_{00}+2 \alpha^{2} R r\right)+\alpha \Pi\left(r_{0}+s_{0}\right)\right\} b^{i}.
\end{aligned}
$$
The constituent functions depending on $\phi = \phi(b^2, s)$ are defined as:
$$
Q := \frac{\phi_{s}}{\phi-s \phi_{s}}, \quad \Theta := \frac{Q-s Q_{s}}{2\left(1+s Q+\left(b^{2}-s^{2}\right) Q_{s}\right)}, \quad \Psi := \frac{Q_{s}}{2\left(1+s Q+\left(b^{2}-s^{2}\right) Q_{s}\right)},
$$
$$
R := \frac{\phi_{b^2}}{\phi-s \phi_{s}}, \quad \Pi := \frac{R_{s}(1+s Q)-2 s R Q_{s}}{1+s Q+\left(b^{2}-s^{2}\right) Q_{s}}, \quad \Omega := \frac{2 R}{1+s Q}-\frac{s+b^{2} Q}{1+s Q} \Pi,
$$
where the subscript $b^2$ denotes the partial derivative with respect to $b^2$.

The explicit  formula of  the Landsberg tensor $L_{j k \ell}$ is given by
\begin{equation}
L_{j k \ell} = -\frac{\rho}{6 \alpha^{5}} \left\{h_{j} h_{k} C_{\ell} + h_{j} h_{\ell} C_{k} + h_{k} h_{\ell} C_{j} + 3 E_{j} h_{k \ell} + 3 E_{k} h_{j \ell} + 3 E_{\ell} h_{j k}\right\}, 
\end{equation}
where the structural blocks and scalars are defined as follows:
$$
\begin{aligned}
& h_{j} := \alpha b_{j}-s y_{j}, \quad h_{j k} := \alpha^{2} a_{j k}-y_{j} y_{k}, \quad \rho := \phi\left(\phi-s \phi_{s}\right), \quad \nu := -Q_{ss}, \\
& \Delta := 1+s Q+\left(b^{2}-s^{2}\right) Q_{s}, \quad \sigma := \frac{1}{\Delta}, \quad \tau := -\frac{Q}{\Delta}, \quad \mu := -\frac{1}{3}\left(Q-s Q_{s}\right), \\
& C_{j} := \left(A_{4} r_{00}+T_{4} \alpha^{2} r-D_{4} \alpha r_{0}-\alpha B_{4} s_{0}\right) h_{j} + 3 \nu J_{j} - \alpha K_{4} M_{j}, \\
& E_{j} := \left(A_{6} r_{00}+T_{6} \alpha^{2} r-D_{6} \alpha r_{0}-\alpha B_{6} s_{0}\right) h_{j} + 3 \mu J_{j} - \alpha K_{6} M_{j}, \\
& J_{j} := \alpha^{2}\left(s_{j 0}+\sigma r_{j 0}+\tau \alpha s_{j}\right)-\left(\sigma r_{00}+\tau \alpha s_{0}\right) y_{j}, \\
& M_{j} := \alpha^{2}\left(r_{j}+s_{j}\right)-\left(r_{0}+s_{0}\right) y_{j}.
\end{aligned}
$$
The auxiliary parameter functions are given by:
$$
\begin{aligned}
A_{6} &= \frac{1}{2 \Delta^{2}}\left\{\left(Q-s Q_{s}\right)^{2}+\left[2\left(s+b^{2} Q\right)-\left(b^{2}-s^{2}\right)\left(Q-s Q_{s}\right)\right] Q_{ss}\right\}, \\
A_{4} &= \frac{1}{2 \Delta^{2}}\left\{-2 \Delta Q_{sss}+3\left(Q-s Q_{s}\right) Q_{ss}+3\left(b^{2}-s^{2}\right) Q_{ss}^{2}\right\}, \\
B_{6} &= D_{6}+2 Q A_{6}-2 Q_{s} \Theta, \quad B_{4} = D_{4}+2 Q A_{4}-6 Q_{ss} \Psi, \\
D_{6} &= \frac{1}{\Delta}\left[\left(b^{2}-s^{2}\right) R_{s}+2 s R\right] s Q_{ss}+\Omega_{s}-s R_{ss}+R_{s}, \\
D_{4} &= -\left\{\left(s+Q b^{2}\right) \Pi_{sss}+3 Q \Omega_{ss}+(1+Q s) \Omega_{sss}-R_{sss}\right\}, \\
K_{6} &= 3 \mu \sigma\left[\left(b^{2}-s^{2}\right) R_{s}+2 s R\right], \quad K_{4} = 3 \nu \sigma\left[\left(b^{2}-s^{2}\right) R_{s}+2 s R\right], \\
T_{6} &= -\frac{1}{\Delta^{2}}\left\{R_{s}\left(Q-s Q_{s}\right) \Delta-\left[\left(b^{2}-s^{2}\right)\left(Q+s Q_{s}\right)+2 s(1+s Q)\right] R Q_{ss}-R\left(Q-s Q_{s}\right)^{2}\right\}, \\
T_{4} &= -\frac{1}{\Delta^{2}}\left\{3 R_{s} Q_{ss} \Delta-3\left[\left(b^{2}-s^{2}\right) Q_{ss}+\left(Q-s Q_{s}\right)\right] R Q_{ss}+2 R Q_{sss} \Delta\right\} .
\end{aligned}
$$

By \cite{Zhou-Wang-Li},  a   reduction lemma for the Landsberg condition  of a general $(\alpha, \beta)$-metric is given as follows:

\begin{lemma}\cite{Zhou-Wang-Li}
Let $F = \alpha \phi\left(b^{2}, s\right)$ be a non-Riemannian general $(\alpha, \beta)$-metric on a manifold of dimension $n \geqslant 3$. Then, the metric is Landsbergian  if and only if the vector components vanish identically:
$$C_{j} = 0 \quad \text{\textit{and}} \quad E_{j} = 0.$$
\end{lemma}
In what follows, we investigate the non-regular, non-Berwaldian Landsberg general $(\alpha,\beta)$-metrics introduced in \cite{Zhou-Wang-Li}. By computing the corresponding integrals in closed form, we derive formulae for the defining function $\phi$ and verify the non-degeneracy condition.

\section{Class I: $(\alpha,\beta)$-metrics} 

~\par

We start with the classical $(\alpha,\beta)$-metrics, a case initially established in \cite[Theorem 1.2 (2)]{Zhou-Wang-Li}. While this case has been explored in the literature, we formalize it here to comprehensively cover all possible special cases of general $(\alpha,\beta)$-metrics mentioned in \cite{Zhou-Wang-Li}.  

\medskip

We can rewrite the fundamental structural relation of this case as follows:
\begin{equation*}
\phi(s) = c_3 \exp \left( \int \frac{c_1 \sqrt{b^2-s^2} + c_2 s}{b^2 + s\left(c_1 \sqrt{b^2-s^2} + c_2 s\right)} \, ds \right),
\end{equation*}
where $\beta$ satisfies
\begin{equation}\label{Class_I_condition_1}
b = b_0, \quad b_{i \mid j} = k\left(b^2 a_{i j} - b_i b_j\right),
\end{equation}
and $b_0, c_1, c_2, c_4$ are constants with $c_1 \neq 0$ and $c_4 > 0$, and $k = k(x)$ is a non-zero scalar function. Under these conditions, the Finsler metric $F$ is a non-Berwaldian Landsberg metric.

\begin{theorem}\label{Thm:Class_I}
Let $F = \alpha \phi(s)$ be an $(\alpha,\beta)$-metric on an $n$-dimensional manifold ($n \geq 3$) such that 
\begin{equation}\label{Eq:Case_2}
\frac{\phi_s}{\phi} = \frac{c_2 s + c_1 \sqrt{b^2-s^2}}{b^2 + c_2 s^2 + c_1 s \sqrt{b^2-s^2}},
\end{equation}
where $c_1 \neq 0$ and $c_2 \neq -1$ are real constants, and $\beta$ satisfies condition \eqref{Class_I_condition_1}. Then, $F$ is a non-regular, non-Berwaldian Landsberg metric, where $\phi(s)$ is given explicitly by
\begin{equation}\label{Eq:Class_2_phi}
\phi(s) = c_3\sqrt{b^2 + c_2 s^2 + c_1 s \sqrt{b^2-s^2}} \, \exp\left( \frac{c_1 }{\sqrt{c_1^2 - 4c_2 - 4}} \operatorname{arctanh} \left( \frac{c_1 s + 2\sqrt{b^2-s^2}}{s \sqrt{c_1^2 - 4c_2 - 4}} \right) \right),
\end{equation}
for some non-zero real constant $c_3$, provided that $c_1^2 - 4c_2 - 4>0$.
\end{theorem}

\begin{proof}
The proof that the metric defined by \eqref{Eq:Case_2} is a non-Berwaldian Landsberg metric was done in \cite[Theorem 1.2 (2)]{Zhou-Wang-Li}. Therefore, our focus here is on deriving an explicit formula for the metric function $\phi$.
Let $\mathcal{D}$ denote the denominator on the right-hand side of \eqref{Eq:Case_2}:
\begin{equation}
\mathcal{D} = b^2 + c_2 s^2 + c_1 s \sqrt{b^2-s^2}.
\end{equation}
Integrating \eqref{Eq:Case_2} with respect to $s$ yields
\begin{equation}\label{Eq:ln_phi_int}
\ln \phi(s) = \int \frac{c_2 s + c_1 \sqrt{b^2-s^2}}{\mathcal{D}} \, ds.
\end{equation}
Differentiating $\mathcal{D}$ with respect to $s$ gives
\begin{equation}
\mathcal{D}_s = 2c_2 s + c_1 \sqrt{b^2-s^2} - \frac{c_1 s^2}{\sqrt{b^2-s^2}} = 2c_2 s + \frac{c_1(b^2 - 2s^2)}{\sqrt{b^2-s^2}}.
\end{equation}
By grouping terms, the original integrand can be decomposed in the following form 
\begin{equation}\label{Eq:Case_2_decomp}
\frac{c_2 s + c_1 \sqrt{b^2-s^2}}{\mathcal{D}} = \frac{1}{2} \frac{\mathcal{D}_s}{\mathcal{D}} + \frac{1}{2} \frac{c_1 b^2}{\sqrt{b^2-s^2} \mathcal{D}}.
\end{equation}
Substituting \eqref{Eq:Case_2_decomp} into \eqref{Eq:ln_phi_int} transforms the integral into
\begin{equation}
\ln \phi(s) = \frac{1}{2} \ln \mathcal{D} + \frac{c_1 b^2}{2} \int \frac{1}{\sqrt{b^2-s^2}\left(b^2 + c_2 s^2 + c_1 s \sqrt{b^2-s^2}\right)} \, ds + \ln c_3.
\end{equation}

To evaluate the remaining transcendental integral, we employ an appropriate trigonometric substitution. Then, we get 
\begin{equation}
\int \frac{b^2}{\sqrt{b^2-s^2}\left(b^2 + c_2 s^2 + c_1 s \sqrt{b^2-s^2}\right)} \, ds = \frac{2}{\sqrt{c_1^2 - 4c_2 - 4}} \operatorname{arctanh} \left( \frac{c_1 s + 2\sqrt{b^2-s^2}}{s \sqrt{c_1^2 - 4c_2 - 4}} \right).
\end{equation}
Multiplying by the coefficient $\frac{c_1 }{2}$ yields the exact exponential argument in \eqref{Eq:Class_2_phi}. Exponentiating both sides completes the derivation of $\phi(s)$.

Substituting \eqref{Eq:Case_2} into the non-degeneracy equation \eqref{Eq:simple_form_almost regularity_1} yields
\begin{equation}
\phi - s \phi_s = \left( 1 - \frac{c_2 s^2 + c_1 s \sqrt{b^2-s^2}}{b^2 + c_2 s^2 + c_1 s \sqrt{b^2-s^2}} \right)\phi = \left( \frac{b^2}{b^2 + c_2 s^2 + c_1 s \sqrt{b^2-s^2}} \right)\phi = \frac{b^2 \phi}{\mathcal{D}}.
\end{equation}
Because $b \neq 0$ and $\phi \neq 0$, we have $\phi - s \phi_s \neq 0$.

For the second non-degeneracy condition \eqref{Eq:simple_form_almost regularity_2}, substituting \eqref{Eq:Case_2} and computing its derivative yields
\begin{equation}
\phi - s \phi_s + (b^2 - s^2)\phi_{ss} = \frac{b^2 (1 + c_2) \phi}{b^2 + c_2 s^2 + c_1 s \sqrt{b^2-s^2}} = \frac{b^2(1+c_2)\phi}{\mathcal{D}}.
\end{equation}
Since $c_2 \neq -1$ by assumption, it follows that $\phi - s \phi_s + (b^2 - s^2)\phi_{ss} \neq 0$.

Consequently, both non-degeneracy conditions hold, which proves that the Finsler metric is non-degenerate.
\end{proof}
Using the fact that 
\begin{equation}
\label{arctan}
\arctan( z)=-\textbf{i} \operatorname{arctanh}(iz), \quad \textbf{i}=\sqrt{-1}
\end{equation}
we have the following corollary.
\begin{corollary}\label{Cor:Class_I_trig}
Under the assumptions of Theorem~\ref{Thm:Class_I}, if the structural constants satisfy the inequality $c_1^2 - 4c_2 - 4 < 0$, then the explicit algebraic form of the metric function $\phi(s)$ transitions into the trigonometric form:
\begin{equation}\label{Eq:Class_2_phi_arctan}
\phi(s) = c_3 \sqrt{b^2 + c_2 s^2 + c_1 s \sqrt{b^2-s^2}} \, \exp\left(- \frac{c_1 }{\sqrt{4c_2 + 4 - c_1^2}} \operatorname{arctan} \left( \frac{c_1 s + 2\sqrt{b^2-s^2}}{s \sqrt{4c_2 + 4 - c_1^2}} \right) \right).
\end{equation}
\end{corollary}

\subsection{Case analysis for $c_1^2 - 4c_2 - 4 = 0$}

When the constants satisfy the algebraic constraint $c_1^2 - 4c_2 - 4 = 0$, the inverse hyperbolic tangent structure from the general case encounters a singular limit. We analyze this special configuration below.

\begin{theorem}
Let $F = \alpha \phi(s)$ be the $(\alpha,\beta)$-metric given in Theorem~\ref{Thm:Class_I}, subject to the algebraic constraint
\begin{equation}\label{eq:constraint_2}
c_1^2 - 4c_2 - 4 = 0, \quad c_1\neq 0 .
\end{equation}
Then, $F$ is a non-regular Landsberg non-Berwaldian metric where the explicit analytical solution $\phi(s)$ reduces to 
\begin{equation}\label{eq:solution_phi_2}
\phi(s) = c_3 \left( c_1 s + 2\sqrt{b^2-s^2} \right) \exp\left( \frac{c_1 s}{c_1 s + 2\sqrt{b^2-s^2}} \right).
\end{equation}
\end{theorem}

\begin{proof}
Substituting the constraint $c_2 = \frac{c_1^2-4}{4}$ into the denominator $\mathcal{D} = b^2 + c_2 s^2 + c_1 s \sqrt{b^2-s^2}$ and rewriting $b^2$ as $(b^2 - s^2) + s^2$, we get
\begin{equation}\label{eq:D_factored_2}
\mathcal{D} = (b^2 - s^2) + c_1 s \sqrt{b^2-s^2} + \frac{c_1^2}{4}s^2 = \left( \frac{c_1}{2}s + \sqrt{b^2-s^2} \right)^2.
\end{equation}
Substituting this form and the value of $c_2$ back into the metric relation yields the following differential equation:
\begin{equation}\label{eq:ode_singular}
\frac{\phi_s}{\phi} = \frac{\frac{c_1^2-4}{4}s + c_1 \sqrt{b^2-s^2}}{\left( \frac{c_1}{2}s + \sqrt{b^2-s^2} \right)^2} = \frac{c_1^2 s-4s+4c_1\sqrt{b^2-s^2} }{\left( c_1s +2 \sqrt{b^2-s^2} \right)^2}.
\end{equation}
This  can be rewritten in the form
\begin{equation}\label{eq:ode_singular_2}
\frac{\phi_s}{\phi} = \frac{c_1 \sqrt{b^2-s^2}-2s}{ (c_1s +2 \sqrt{b^2-s^2}) \sqrt{b^2-s^2}}+\frac{2c_1b^2}{\left( c_1s +2 \sqrt{b^2-s^2} \right)^2\sqrt{b^2-s^2}}.
\end{equation}
The first fraction can be seen as $\frac{\mathcal{D}_s}{\mathcal{D}}$. For the second fraction, one can see that 
$$u:= \frac{c_1 s}{c_1 s + 2\sqrt{b^2-s^2}}, \quad u_s=\frac{2c_1b^2}{\left( c_1s +2 \sqrt{b^2-s^2} \right)^2\sqrt{b^2-s^2}}.$$
Integrating both sides yields the desired expression for  $\phi$. Regarding the non-degeneracy property, this can be verified via Theorem~\ref{Thm:Class_I}.
\end{proof}

\section{Class II}
~\par
Following the classification established in \cite[Theorem 1.2 (3)]{Zhou-Wang-Li}, there exists a function $\tau = \exp \left( \int \frac{k-c}{2 b^2 c} \, db^2 \right)$ transforming the structural quantities into invariant forms: $\widetilde{\beta} = \tau \beta$, $\widetilde{b}^2 = \tau^2 b^2$, and $\widetilde{s} = \tau s$. Under this framework, the metric function $\phi(b^2, s)$ can be expressed solely in terms of $\widetilde{b}^2$ and $\widetilde{s}$ as:
\begin{equation*}
\phi\left(\widetilde{b}^2, \widetilde{s}\right)=\lambda_3 \exp \left( \int  \frac{2 c_1 \sqrt{\widetilde{b}^2-\widetilde{s}^2}+\widetilde{b}^2 \widetilde{s} \lambda_1+2 \widetilde{b}^2 \widetilde{s} \lambda_0}{2 \widetilde{s} c_1 \sqrt{\widetilde{b}^2-\widetilde{s}^2}+\widetilde{b}^2 \widetilde{s}^2 \lambda_1+\widetilde{b}^2-2 \widetilde{b}^2\left(\widetilde{b}^2-\widetilde{s}^2\right) \lambda_0} \, d\widetilde{s} \right),
\end{equation*}
where $\widetilde{\beta}=\widetilde{b}_i(x) y^i$ satisfies 
$$\widetilde{b}_{i \mid j}=\widetilde{c} a_{i j},$$ 
with $\widetilde{c}=b^2 k \exp \left( \int \frac{k-c}{2 b^2 c} \, db^2 \right)$. The parameters $c_1 = c_1(\widetilde{b}^2)$, $\lambda_0 = \lambda_0(\widetilde{b}^2)$, $\lambda_1 = \lambda_1(\widetilde{b}^2)$, and $\lambda_3 = \lambda_3(\widetilde{b}^2) > 0$ are smooth functions satisfying the standard derivative compatibility conditions:
\begin{align}
\label{Class_II_condition}
\nonumber\int  A_1\left(\widetilde{b}^2, \widetilde{s}\right) d \widetilde{s}= & \frac{1}{2 \widetilde{b}^4 \widetilde{s} \lambda_3}\left\{\left[\widetilde{b}^4\left(\widetilde{b}^2-\widetilde{s}^2\right)\left(2 \lambda_0+\widetilde{s}^2 \lambda_2\right)-2 \widetilde{s} c_1\left(\widetilde{b}^2-\widetilde{s}^2\right)^{\frac{3}{2}}-\widetilde{b}^4\right] \lambda_3 A\left(\widetilde{b}^2, \widetilde{s}\right)\right. \\
& \left.+\left[\widetilde{b}^4 \widetilde{s}\left(\lambda_2 \widetilde{s}^2+2 \lambda_1+2 \lambda_0\right)+2 c_1\left(\widetilde{b}^2-\widetilde{s}^2\right)^{\frac{3}{2}}\right] \lambda_3-2 \widetilde{b}^4 \widetilde{s} \lambda_3^{\prime}\right\}, 
\end{align}
where
\begin{align*}
&A\left(\widetilde{b}^2, \widetilde{s}\right)=\frac{2 c_1 \sqrt{\widetilde{b}^2-\widetilde{s}^2}+\widetilde{b}^2 \widetilde{s} \lambda_1+2 \widetilde{b}^2 \widetilde{s} \lambda_0}{2 \widetilde{s} c_1 \sqrt{\widetilde{b}^2-\widetilde{s}^2}+\widetilde{b}^2 \widetilde{s}^2 \lambda_1+\widetilde{b}^2-2 \widetilde{b}^2\left(\widetilde{b}^2-\widetilde{s}^2\right) \lambda_0}, \quad A_1\left(\widetilde{b}^2, \widetilde{s}\right)=\frac{\partial A\left(\widetilde{b}^2, \widetilde{s}\right)}{\partial \widetilde{b}^2}.
\end{align*}

\begin{theorem}\label{Thm:Class_II}
Let $F=\alpha \phi(\tilde{b}^2, \tilde{s})$ be a general $(\alpha,\beta)$-metric on an $n$-dimensional manifold ($n \geq 3$) such that 
\begin{equation}\label{Eq:Case_1}
\frac{\phi_{\tilde{s}}}{\phi} = \frac{2 c_1 \sqrt{\tilde{b}^2-\tilde{s}^2} + \tilde{b}^2 \tilde{s} \lambda_1 + 2 \tilde{b}^2 \tilde{s} \lambda_0}{2 \tilde{s} c_1 \sqrt{\tilde{b}^2-\tilde{s}^2} + \tilde{b}^2 \tilde{s}^2 \lambda_1 + \tilde{b}^2 - 2 \tilde{b}^2(\tilde{b}^2-\tilde{s}^2)\lambda_0}, 
\end{equation}
where $c_1 \in \mathbb{R}$ is a constant, and $\lambda_0 = \lambda_0(\tilde{b}^2)$, $\lambda_1 = \lambda_1(\tilde{b}^2)$ depend only on $\tilde{b}^2$ with condition \eqref{Class_II_condition} satisfied, as well as $1-2\tilde{b}^2\lambda_0\neq 0$ and $\tilde{b}^2\lambda_1+1\neq 0$. Then, $F$ is a non-regular Landsberg non-Berwaldian metric. Moreover, $\phi$ is explicitly given by 
\begin{equation}\label{Eq:Class_1_phi}
\phi(\tilde{b}^2, \tilde{s}) = \lambda_3(\tilde{b}^2) \sqrt{2 \tilde{s} c_1 \sqrt{\tilde{b}^2-\tilde{s}^2} + \tilde{b}^2 \tilde{s}^2 \lambda_1 + \tilde{b}^2 - 2\tilde{b}^2(\tilde{b}^2-\tilde{s}^2)\lambda_0}\ e^{-\frac{c_1}{\sqrt{\Delta_2}} \operatorname{arctanh} \left( \frac{(\tilde{b}^2\lambda_1+1)\tilde{s} + c_1\sqrt{\tilde{b}^2-\tilde{s}^2}}{\sqrt{\Delta_2}\sqrt{\tilde{b}^2-\tilde{s}^2}} \right) },
\end{equation}
where $\lambda_3(\tilde{b}^2)$ is a non-zero smooth function, and $\Delta_2$ is defined by
\begin{equation}\label{Eq:Class_1_3}
\Delta_2 =  c_1^2 -(\tilde{b}^2\lambda_1 + 1)(1 - 2\tilde{b}^2\lambda_0)>0.
\end{equation}
\end{theorem}

\begin{proof}
The proof that the metric determined by \eqref{Eq:Case_1} is a non-Berwaldian Landsberg metric was shown in \cite[Theorem 1.2 (3)]{Zhou-Wang-Li}; therefore, this analysis focuses on deriving an explicit formula for the metric function \(\phi\).
Let $\mathcal{D}$ denote the denominator of the right-hand side of \eqref{Eq:Case_1}:
\begin{equation}
\mathcal{D} = 2 \tilde{s} c_1 \sqrt{\tilde{b}^2-\tilde{s}^2} + \tilde{b}^2 \tilde{s}^2 \lambda_1 + \tilde{b}^2 - 2 \tilde{b}^2(\tilde{b}^2-\tilde{s}^2)\lambda_0.
\end{equation}
Integrating \eqref{Eq:Case_1} with respect to $\tilde{s}$ yields 
\begin{equation}\label{Eq:ln_phi_int_I}
\ln \phi = \int \frac{2 c_1 \sqrt{\tilde{b}^2-\tilde{s}^2} + \tilde{b}^2 \tilde{s} \lambda_1 + 2 \tilde{b}^2 \tilde{s} \lambda_0}{\mathcal{D}} \, d\tilde{s}.
\end{equation}
Differentiating $\mathcal{D}$ directly with respect to $\tilde{s}$ yields:
\begin{equation}
\mathcal{D}_s = 2c_1\sqrt{\tilde{b}^2-\tilde{s}^2} - \frac{2c_1\tilde{s}^2}{\sqrt{\tilde{b}^2-\tilde{s}^2}} + 2\tilde{b}^2\tilde{s}\lambda_1 + 4\tilde{b}^2\tilde{s}\lambda_0.
\end{equation}
Using the algebraic identity $\tilde{s}^2 = \tilde{b}^2 - (\tilde{b}^2 - \tilde{s}^2)$, we decompose the integrand into the following decomposition
\begin{equation}\label{Eq:Case_1_decomp}
\frac{2 c_1 \sqrt{\tilde{b}^2-\tilde{s}^2} + \tilde{b}^2 \tilde{s} \lambda_1 + 2 \tilde{b}^2 \tilde{s} \lambda_0}{\mathcal{D}} = \frac{1}{2} \frac{\mathcal{D}_s}{\mathcal{D}} + \frac{\tilde{b}^2 c_1}{\sqrt{\tilde{b}^2-\tilde{s}^2} \mathcal{D}}.
\end{equation}
Substituting \eqref{Eq:Case_1_decomp} back into \eqref{Eq:ln_phi_int_I} gives:
\begin{equation}
\ln \phi = \frac{1}{2} \ln |\mathcal{D}| + \tilde{b}^2 c_1 \int \frac{1}{\sqrt{\tilde{b}^2-\tilde{s}^2} \mathcal{D}} \, d\tilde{s} + \ln \lambda_3(\tilde{b}^2).
\end{equation}
To evaluate the remaining integral, we introduce the substitution variable $u = \frac{\tilde{s}}{\sqrt{\tilde{b}^2-\tilde{s}^2}}$, which yields 
$$d\tilde{s} = \frac{\tilde{b}^2}{(\tilde{b}^2(1+u^2))^{3/2}}\,du = \frac{1}{\tilde{b}(1+u^2)^{3/2}}\,du,$$
 simplifying the expression to:
\begin{equation}
\int \frac{1}{\sqrt{\tilde{b}^2-\tilde{s}^2} \mathcal{D}} \, d\tilde{s} = \int \frac{1}{\tilde{b}^2(1-2\tilde{b}^2\lambda_0) + 2\tilde{b}^2 c_1 u + \tilde{b}^2(\tilde{b}^2\lambda_1 + 1)u^2} \, du.
\end{equation}
Factoring out $\tilde b^2$ from the denominator gives
\begin{equation}
c_1 \int \frac{1}{(\tilde{b}^2\lambda_1 + 1)u^2 + 2c_1 u + (1-2\tilde{b}^2\lambda_0)} \, du = -\frac{c_1}{\sqrt{\Delta_2}}\operatorname{arctanh}\left(\frac{(\tilde{b}^2\lambda_1+1)u + c_1}{\sqrt{\Delta_2}}\right).
\end{equation}
Substituting $u = \frac{\tilde{s}}{\sqrt{\tilde{b}^2-\tilde{s}^2}}$ back into the expression produces the exact functional components within \eqref{Eq:Class_1_phi}. Exponentiating both sides completes the proof.

To ensure the geometric non-degeneracy of the resulting Landsberg metric, we evaluate the fundamental non-degeneracy conditions. Substituting \eqref{Eq:Case_1} directly into the non-degeneracy condition \eqref{Eq:simple_form_almost regularity_1} yields:
\begin{equation}
\phi - \tilde{s}\phi_{\tilde{s}} =\left( \frac{\mathcal{D} - \tilde{s}\left(2 c_1 \sqrt{\tilde{b}^2-\tilde{s}^2} + \tilde{b}^2 \tilde{s} \lambda_1 + 2 \tilde{b}^2 \tilde{s} \lambda_0\right)}{\mathcal{D}} \right)  \phi = \frac{\tilde{b}^2(1-2\tilde{b}^2\lambda_0)\phi}{\mathcal{D}}.
\end{equation}
Because $\tilde{b} > 0$, $\phi > 0$, and $1 - 2\tilde{b}^2\lambda_0 \neq 0$, the metric satisfies $\phi - \tilde{s}\phi_{\tilde{s}} \neq 0$.

Next, evaluating the second non-degeneracy condition \eqref{Eq:simple_form_almost regularity_2}   yields:
\begin{equation}
\phi - \tilde{s}\phi_{\tilde{s}} + (\tilde{b}^2 - \tilde{s}^2)\phi_{\tilde{s}\tilde{s}} = \frac{\tilde{b}^2(1-2\tilde{b}^2\lambda_0)(\tilde{b}^2\lambda_1+1)\phi}{\mathcal{D}^2}.
\end{equation}
Given that $1-2\tilde{b}^2\lambda_0 \neq 0$ and $\tilde{b}^2\lambda_1+1 \neq 0$, then $\phi - \tilde{s}\phi_{\tilde{s}} + (\tilde{b}^2 - \tilde{s}^2)\phi_{\tilde{s}\tilde{s}}\neq0$. This mathematically confirms that the general solution represents a  non-degenerate, non-regular Finsler metric.
\end{proof}
By using \eqref{arctan}, we get the following formula.
\begin{corollary}\label{Cor:Class_II_trig}
Under the assumptions of Theorem~\ref{Thm:Class_II}, if   $\Delta_2 < 0$, then the explicit algebraic form of the metric function $\phi(\tilde{b}^2, \tilde{s})$ transitions into the  form  
\begin{equation}\label{Eq:Class_1_phi_arctan}
\phi(\tilde{b}^2, \tilde{s}) = \lambda_3(\tilde{b}^2) \sqrt{2 \tilde{s} c_1 \sqrt{\tilde{b}^2-\tilde{s}^2} + \tilde{b}^2 \tilde{s}^2 \lambda_1 + \tilde{b}^2 - 2\tilde{b}^2(\tilde{b}^2-\tilde{s}^2)\lambda_0}\ e^{\frac{c_1}{\sqrt{-\Delta_2}} \operatorname{arctan} \left( \frac{(\tilde{b}^2\lambda_1+1)\tilde{s} + c_1\sqrt{\tilde{b}^2-\tilde{s}^2}}{\sqrt{-\Delta_2}\sqrt{\tilde{b}^2-\tilde{s}^2}} \right) }.
\end{equation}
\end{corollary}

\subsection{Case Analysis for $\Delta_2 = 0$}
\begin{theorem}
Let $F=\alpha\phi(\tilde{b}^2, \tilde{s})$ be the general $(\alpha,\beta)$-metric given in Theorem~\ref{Thm:Class_II} subject to the algebraic   constraint:
\begin{equation}\label{eq:constraint}
\Delta_2 = c_1^2 - (\tilde{b}^2\lambda_1 + 1)(1 - 2\tilde{b}^2\lambda_0) = 0.
\end{equation}
Then, the explicit analytical solution of $\phi$ reduces to 
\begin{equation}\label{eq:solution_phi}
\phi(\tilde{b}^2, \tilde{s}) = \lambda_3(\tilde{b}^2)\left( (\tilde{b}^2\lambda_1+1)\tilde{s} + c_1\sqrt{\tilde{b}^2-\tilde{s}^2} \right) \exp\left( - \frac{c_1 \sqrt{\tilde{b}^2-\tilde{s}^2}}{(\tilde{b}^2\lambda_1+1)\tilde{s} + c_1\sqrt{\tilde{b}^2-\tilde{s}^2}} \right).
\end{equation}
\end{theorem}

\begin{proof}
Let $\mathcal{D}$ denote the denominator of \eqref{Eq:Case_1}. Grouping terms by scaling powers of $\tilde{s}^2$ yields:
\begin{equation}\label{eq:D_unfactored}
\mathcal{D} = 2 c_1 \tilde{s} \sqrt{\tilde{b}^2-\tilde{s}^2} + \tilde{s}^2(\tilde{b}^2\lambda_1 + 2\tilde{b}^2\lambda_0) + \tilde{b}^2(1 - 2\tilde{b}^2\lambda_0).
\end{equation}
Isolating the boundary factor from the condition $\Delta_2 = 0$ given in \eqref{eq:constraint} shows that:
\begin{equation}\label{eq:constraint_sub}
1 - 2\tilde{b}^2\lambda_0 = \frac{c_1^2}{\tilde{b}^2\lambda_1 + 1} \implies \tilde{b}^2(1 - 2\tilde{b}^2\lambda_0) = \frac{\tilde{b}^2 c_1^2}{\tilde{b}^2\lambda_1 + 1}.
\end{equation}
Substituting \eqref{eq:constraint_sub} into \eqref{eq:D_unfactored} and using the identity $\tilde{s}^2(\tilde{b}^2\lambda_1+2\tilde{b}^2\lambda_0) = \tilde{s}^2(\tilde{b}^2\lambda_1+1) - \tilde{s}^2(1-2\tilde{b}^2\lambda_0)$, we observe that the denominator collects into a perfect-square algebraic factorization:
\begin{equation}\label{eq:D_factored}
\mathcal{D} = \frac{1}{\tilde{b}^2\lambda_1 + 1} \left( (\tilde{b}^2\lambda_1 + 1)\tilde{s} + c_1\sqrt{\tilde{b}^2-\tilde{s}^2} \right)^2.
\end{equation}
By evaluating the primary differential equation under this specific constraint, the structural ratio $\frac{\phi_{\tilde{s}}}{\phi}$ simplifies to:
\begin{equation}\label{eq:ratio_special_case}
\frac{\phi_{\tilde{s}}}{\phi} = \frac{(\tilde{b}^2\lambda_1 + 1)\left(2 c_1 \sqrt{\tilde{b}^2-\tilde{s}^2} + \tilde{b}^2 \tilde{s} \lambda_1 + 2 \tilde{b}^2 \tilde{s} \lambda_0\right)}{\left( (\tilde{b}^2\lambda_1 + 1)\tilde{s} + c_1\sqrt{\tilde{b}^2-\tilde{s}^2} \right)^2}.
\end{equation}
To perform the integration of this ratio with respect to $\tilde{s}$, we reuse the standard coordinate transformation $u = \frac{\tilde{s}}{\sqrt{\tilde{b}^2-\tilde{s}^2}}$, turning the integration problem into:
\begin{equation}
\ln \phi = \int \frac{\phi_{\tilde{s}}}{\phi} \, d\tilde{s} = \frac{1}{2}\ln\mathcal{D} + \int \frac{c_1}{(\tilde{b}^2\lambda_1+1)u^2 + 2c_1 u + (1-2\tilde{b}^2\lambda_0)}\,du.
\end{equation}
Because $\Delta_2 = 0$, the quadratic expression in the denominator under the integral sign is a perfect square, given by $\frac{1}{\tilde{b}^2\lambda_1+1}\left((\tilde{b}^2\lambda_1+1)u + c_1\right)^2$. Thus, the integration simplifies directly to:
\begin{equation}
c_1(\tilde{b}^2\lambda_1+1)\int \frac{1}{\left((\tilde{b}^2\lambda_1+1)u+c_1\right)^2}\,du = -\frac{c_1}{\left((\tilde{b}^2\lambda_1+1)u + c_1\right)}.
\end{equation}
Substituting $u = \frac{\tilde{s}}{\sqrt{\tilde{b}^2-\tilde{s}^2}}$ back into this integrated term yields:
\begin{equation}
-\frac{c_1}{\left((\tilde{b}^2\lambda_1+1)u + c_1\right)}=\frac{-c_1 \sqrt{\tilde{b}^2-\tilde{s}^2}}{(\tilde{b}^2\lambda_1+1)\tilde{s} + c_1\sqrt{\tilde{b}^2-\tilde{s}^2}}.
\end{equation}
Combining this resulting function with $\frac{1}{2}\ln\mathcal{D}$ (where the constant scale factor $1/\sqrt{\tilde{b}^2\lambda_1+1}$ is absorbed into the arbitrary function $\ln \lambda_3(\tilde{b}^2)$), we obtain:
\begin{equation}\label{eq:ln_phi}
\ln \phi = \ln \left( (\tilde{b}^2\lambda_1+1)\tilde{s} + c_1\sqrt{\tilde{b}^2-\tilde{s}^2} \right) - \frac{c_1 \sqrt{\tilde{b}^2-\tilde{s}^2}}{(\tilde{b}^2\lambda_1+1)\tilde{s} + c_1\sqrt{\tilde{b}^2-\tilde{s}^2}} + \ln \lambda_3(\tilde{b}^2).
\end{equation}
Using the algebraic decomposition $\frac{c_1\sqrt{\tilde{b}^2-\tilde{s}^2}}{(\tilde{b}^2\lambda_1+1)\tilde{s} + c_1\sqrt{\tilde{b}^2-\tilde{s}^2}} = 1 - \frac{(\tilde{b}^2\lambda_1+1)\tilde{s}}{(\tilde{b}^2\lambda_1+1)\tilde{s} + c_1\sqrt{\tilde{b}^2-\tilde{s}^2}}$, the constant $1$ is absorbed into the arbitrary smooth function $\ln \lambda_3(\tilde{b}^2)$. Exponentiating both sides matches the explicit form of \eqref{eq:solution_phi}.

The non-degeneracy property is already verified in the proof of Theorem~\ref{Thm:Class_II}.
\end{proof}

\section{Classes III and IV}
~\par

For further analysis and structural study of the third and fourth non-Berwaldian Landsberg general $(\alpha,\beta)$-metrics, it is necessary to establish an explicit relationship between the functional coefficients $f_1(b^2)$ and $f_2(b^2)$. This underlying relation is formalized by the following lemma.

\begin{lemma}
Let $f_1(b^2)$ and $f_2(b^2)$ be smooth functional terms defined as follows:
\begin{align}
f_1(b^2) &= c_1 \exp\left( \int -\frac{1 + \lambda b^2 + b^6 f(b^2)}{b^2(1 + \lambda b^2)} \, db^2 \right), \label{eq:f1_def} \\
f_2(b^2) &= c_2 \exp\left( \int -\frac{2 b^4 f(b^2)}{1 + \lambda b^2} \, db^2\right), \label{eq:f2_def}
\end{align}
where $c_1, c_2$ are arbitrary integration constants, and $\lambda$ is an arbitrary function of $b^2$. Then, these functions satisfy the explicit  algebraic relationship
\begin{equation}\label{eq:f1_f2_relation}
f_2 = C b^4 f_1^2,
\end{equation}
where $C = \frac{c_2}{c_1^2}$.
\end{lemma}

\begin{proof}
We begin by decomposing the integrand of $f_1(b^2)$ in \eqref{eq:f1_def} as follows:
\begin{equation}\label{eq:integrand_decomp}
\frac{1 + \lambda b^2 + b^6 f(b^2)}{b^2(1 + \lambda b^2)} = \frac{1 + \lambda b^2}{b^2(1 + \lambda b^2)} + \frac{b^6 f(b^2)}{b^2(1 + \lambda b^2)} = \frac{1}{b^2} + \frac{b^4 f(b^2)}{1 + \lambda b^2}.
\end{equation}
Substituting \eqref{eq:integrand_decomp} back into the exponential expression for $f_1(b^2)$, we have  
\begin{align}
f_1(b^2) &= c_1 \exp\left( \int -\frac{1}{b^2} \, db^2 + \int -\frac{b^4 f(b^2)}{1 + \lambda b^2} \, db^2 \right) \nonumber \\
&= c_1 \exp(-\ln b^2) \cdot \exp\left( \int -\frac{b^4 f(b^2)}{1 + \lambda b^2} \, db^2 \right) \nonumber \\
&= \frac{c_1}{b^2} \exp\left( \int -\frac{b^4 f(b^2)}{1 + \lambda b^2} \, db^2 \right). \label{eq:f1_isolated}
\end{align}
Rearranging   \eqref{eq:f1_isolated} yields  
\begin{equation}\label{eq:kernel_expression}
\exp\left( \int -\frac{b^4 f(b^2)}{1 + \lambda b^2} \, db^2 \right) = \frac{b^2 f_1(b^2)}{c_1}.
\end{equation}
Next, we examine the structural form of $f_2(b^2)$ given in \eqref{eq:f2_def}. By utilizing the standard properties of the exponential mapping, we have 
\begin{equation}\label{eq:f2_factored}
f_2(b^2) = c_2 \exp\left( 2 \int -\frac{b^4 f(b^2)}{1 + \lambda b^2} \, db^2 \right) = c_2 \left[ \exp\left( \int -\frac{b^4 f(b^2)}{1 + \lambda b^2} \, db^2 \right) \right]^2.
\end{equation}
Direct substitution  of \eqref{eq:kernel_expression} into  \eqref{eq:f2_factored}  implies
\begin{equation}\label{eq:algebraic_b2}
f_2(b^2) = c_2 \left( \frac{b^2 f_1(b^2)}{c_1} \right)^2 = \frac{c_2}{c_1^2} b^4 f_1(b^2)^2.
\end{equation}
Setting   $C = \frac{c_2}{c_1^2}$ simplifies \eqref{eq:algebraic_b2} into  the required formula.
This completes the proof.
\end{proof}

The third class considered here corresponds to the family obtained in \cite[Theorem 1.3 (2)]{Zhou-Wang-Li}, where the structural function $\phi$ is given implicitly by
\begin{equation}\label{phi_Class_3_1}
\phi=f_3\left(b^2\right) \exp \left(\int \frac{b^2 \sqrt{b^2-s^2} f_1\left(b^2\right)+s\left(f_2\left(b^2\right)-1\right)}{b^2+b^2 s \sqrt{b^2-s^2} f_1\left(b^2\right)+s^2\left(f_2\left(b^2\right)-1\right)} \, d s\right),
\end{equation}
and the underlying differential forms of the metric components $\beta=b_i(x) y^i$ satisfy the tensor fields:
\begin{align}\label{Case_III_conditions}
& s_{i j}=\frac{1}{b^2}\left\{b_i s_j-b_j s_i\right\}, \nonumber \\
& r_{i j}=k\left(b^2 a_{i j}-b_i b_j\right)+c b_i b_j+\lambda\left\{b_i s_j+b_j s_i\right\},
\end{align}
where $k=k(x)$, $c=c(x)$, and $\lambda=\lambda(x) \neq -\frac{1}{b^2}$ are scalar functions, $f_1\left(b^2\right)=c_1 \mathrm{e}^{\int-\frac{1+\lambda b^2+b^6 f\left(b^2\right)}{b^2\left(1+\lambda b^2\right)} \, d b^2} \neq 0$, $f_2\left(b^2\right)=c_2 \mathrm{e}^{\int-\frac{2 b^4 f\left(b^2\right)}{\left(1+\lambda b^2\right)} \, d b^2}$, while $f_3\left(b^2\right)$ and $f\left(b^2\right)$ are $C^{\infty}$ smooth functions of $b^2$, and $c_1, c_2$ are real constants.
 
The fourth class treated in this work is the counterpart family obtained in \cite[Theorem 1.4 (2)]{Zhou-Wang-Li}, where $\phi$ maintains the integral representation \eqref{phi_Class_3_1} and the covariant derivative of $\beta$ simplifies to:
\begin{align}\label{Case_VI_conditions} 
b_{i \mid j}=k\left(b^2 a_{i j}-b_i b_j\right)+c b_i b_j+\theta_i b_j+\theta_j b_i,
\end{align}
where $\theta_i$ is a 1-form perpendicular to $\beta$.

\begin{theorem}\label{Thm:Class_III_IV}
Let $F=\alpha \phi(\tilde{b}^2, \tilde{s})$ be a general $(\alpha,\beta)$-metric whose structural function satisfies the first-order ordinary differential equation:
\begin{equation}\label{Eq:Case_3_DE}
\frac{\phi_s}{\phi} = \frac{b^2 \sqrt{b^2-s^2} f_1(b^2) + s(f_2(b^2) - 1)}{b^2 + b^2 s \sqrt{b^2-s^2} f_1(b^2) + s^2(f_2(b^2) - 1)},
\end{equation}
under the  conditions \eqref{Case_III_conditions} or \eqref{Case_VI_conditions}. Then, $F$ is a non-regular Landsberg non-Berwaldian metric if and only if $f_1(b^2) \neq 0$ (or, equivalently, $f_2(b^2) \neq 0$). Moreover, the analytical metric function $\phi$ is explicitly given by:
\begin{equation}\label{Eq:Class_3_phi}
\phi(b^2, s) = f_3(b^2) \sqrt{b^2 + b^2 s \sqrt{b^2-s^2} f_1(b^2) + s^2(f_2(b^2) - 1)} \, e^{ -\frac{b^2 f_1(b^2)}{ \sqrt{\Delta_3}} \operatorname{arctanh} \left( \frac{b^2 f_1(b^2) \sqrt{b^2-s^2} + 2s f_2(b^2)}{\sqrt{\Delta_3} \sqrt{b^2-s^2}} \right) },
\end{equation}
where $f_3(b^2) > 0$ is an arbitrary smooth  function, and $\Delta_3$ is  defined as:
\begin{equation}\label{Eq:Class_3_Omega}
\Delta_3 = b^4 f_1(b^2)^2 - 4f_2(b^2)>0.
\end{equation}
\end{theorem}

\begin{proof}
The proof that the metric provided by \eqref{Eq:Case_3_DE} is a non-Berwaldian Landsberg metric was established in \cite[Theorem 1.3 (2)]{Zhou-Wang-Li} and \cite[Theorem 1.4 (2)]{Zhou-Wang-Li}. Therefore, this analysis focuses on deriving an explicit formula for the metric function $\phi$.
To streamline the algebraic manipulations, we introduce the following
\begin{equation*}
A := f_1(b^2), \quad B := f_2(b^2) - 1.
\end{equation*}
Let $\mathcal{D}$ represent the denominator of \eqref{Eq:Case_3_DE}, hence we have
\begin{equation}\label{eq:Denom_D}
\mathcal{D} = b^2 + b^2 s A \sqrt{b^2-s^2} + s^2 B.
\end{equation}
Rewriting \eqref{Eq:Case_3_DE} as follows
\begin{equation}\label{eq:pde_rewritten}
\frac{\phi_s}{\phi} = \frac{b^2 A \sqrt{b^2-s^2} + s B}{\mathcal{D}}.
\end{equation}
Integrating \eqref{eq:pde_rewritten} with respect to the variable $s$ establishes
\begin{equation}\label{eq:ln_phi_int}
\ln \phi = \int \frac{b^2 A \sqrt{b^2-s^2} + s B}{\mathcal{D}} \, ds + \ln f_3(b^2).
\end{equation}
To resolve this integral structure, we compute the derivative of $\mathcal{D}$ with respect to $s$. We obtain
\begin{equation}
\mathcal{D}_s = \frac{b^2 A(b^2-2s^2)}{\sqrt{b^2-s^2}} + 2s B.
\end{equation}
The numerator of the integrand admits the decomposition
\begin{equation}\label{eq:decomp_main}
\frac{b^2 A \sqrt{b^2-s^2} + s B}{\mathcal{D}} = \frac{1}{2} \frac{\mathcal{D}_s}{\mathcal{D}} + \mathcal{R}.
\end{equation}
By direct evaluation,   $\mathcal{R}$ isolates as follows
\begin{align*}
\mathcal{R} &= \frac{2\left(b^2 A \sqrt{b^2-s^2} + s B\right) - \mathcal{D}_s}{2\mathcal{D}}= \frac{b^2 A}{2 \sqrt{b^2-s^2} \mathcal{D}}.
\end{align*}
Substituting $\mathcal{R}$ back into \eqref{eq:ln_phi_int} yields
\begin{equation}\label{eq:ln_phi_split}
\ln \phi = \frac{1}{2} \ln \mathcal{D} + \frac{b^2 A}{2} \int \frac{1}{\sqrt{b^2-s^2} \mathcal{D}} \, ds + \ln f_3(b^2).
\end{equation}
We identify the secondary integral expression as $I = \int \frac{1}{\sqrt{b^2-s^2} \mathcal{D}} \, ds$. Using the  fact that
\begin{equation}
\frac{d}{ds} \left[ \operatorname{arctanh} \left( \frac{b^2 A \sqrt{b^2-s^2} + 2s f_2(b^2)}{\sqrt{\Delta_3} \sqrt{b^2-s^2}} \right) \right] = -\frac{\sqrt{\Delta_3}}{2 \sqrt{b^2-s^2} \mathcal{D}}.
\end{equation}
Thus, the integral evaluates explicitly to
\begin{equation}\label{eq:integral_eval}
I = -\frac{2}{\sqrt{\Delta_3}} \operatorname{arctanh} \left( \frac{b^2 A \sqrt{b^2-s^2} + 2s f_2(b^2)}{\sqrt{\Delta_3} \sqrt{b^2-s^2}} \right).
\end{equation}
Plugging \eqref{eq:integral_eval} back into the split logarithmic equation \eqref{eq:ln_phi_split}, the front coefficient reduces via $\frac{b^2 A}{2} \left(-\frac{2}{\sqrt{\Delta_3}}\right) = -\frac{b^2 A}{\sqrt{\Delta_3}}$. Exponentiating both sides resolves into the analytical form:
\begin{equation*}
\phi = f_3(b^2) \sqrt{\mathcal{D}} \, \exp \left( -\frac{b^2 A}{\sqrt{\Delta_3}} \operatorname{arctanh} \left( \frac{b^2 A \sqrt{b^2-s^2} + 2s f_2(b^2)}{\sqrt{\Delta_3} \sqrt{b^2-s^2}} \right) \right).
\end{equation*}
Restoring the functional definitions $A = f_1(b^2)$ and $B = f_2(b^2)-1$ yields the explicit metric solution \eqref{Eq:Class_3_phi}.

To verify the geometric validity and non-degeneracy of the obtained metric, we test the core conditions for  non-degeneracy. First, we analyze the first non-degeneracy condition. Substituting the structural derivative ratio \eqref{Eq:Case_3_DE} directly into the non-degeneracy condition \eqref{Eq:simple_form_almost regularity_1} yields 
\begin{equation}
\phi - s\phi_s = \phi \left( \frac{\mathcal{D} - s\left(b^2 A \sqrt{b^2-s^2} + s B\right)}{\mathcal{D}} \right)=\frac{b^2 \phi}{\mathcal{D}}.
\end{equation}
Since $b > 0$, $\phi > 0$, then  the condition $\phi - s\phi_s \neq 0$ is  satisfied.

Next, we evaluate the second non-degeneracy condition \eqref{Eq:simple_form_almost regularity_2}, which give 
\begin{equation}
\phi - s\phi_s + (b^2 - s^2)\phi_{ss} =\frac{b^4 f_2(b^2) \phi}{\mathcal{D}^2} \neq 0.
\end{equation}
Since   $f_2(b^2) \neq 0$, then this   confirms that the general solution represents a  non-degenerate, non-regular  Finsler metric.
\end{proof}
By making use of \eqref{arctan}, we provide the following special case. 
\begin{corollary}\label{Cor:Class_III_IV_trig}
Under the assumptions of Theorem \ref{Thm:Class_III_IV}, if  $\Delta_3 < 0$, then the explicit algebraic form of the analytical metric function $\phi(b^2, s)$ transitions into the following  trigonometric form 
\begin{equation}\label{Eq:Class_3_phi_arctan}
\phi(b^2, s) = f_3 \sqrt{b^2 + b^2 s \sqrt{b^2-s^2} f_1(b^2) + s^2(f_2(b^2) - 1)} \, e^{ \frac{b^2 f_1(b^2)}{ \sqrt{-\Delta_3}} \operatorname{arctan} \left( \frac{b^2 f_1(b^2) \sqrt{b^2-s^2} + 2s f_2(b^2)}{\sqrt{-\Delta_3} \sqrt{b^2-s^2}} \right) }.
\end{equation}
\end{corollary}

\subsection{The singular case $\Delta_3=0$}

\begin{theorem}
Let $F=\alpha\phi(\tilde{b}^2, \tilde{s})$ be the general $(\alpha,\beta)$-metric given in Theorem \ref{Thm:Class_III_IV} subject to the specific algebraic  constraint:
\begin{equation}\label{eq:f2_constraint_sing}
f_2(b^2) = \frac{1}{4} b^4 f_1^2(b^2).
\end{equation}
Then, the explicit analytical solution to  \eqref{Eq:Case_3_DE} is given by:
\begin{equation}\label{eq:explicit_phi_sol_sing}
\phi(b^2, s) = f_3(b^2) \left(  b^2 f_1(b^2) s + 2 \sqrt{b^2-s^2} \right) \exp\left( -\frac{2\sqrt{b^2-s^2} }{b^2 f_1(b^2) s + 2\sqrt{b^2-s^2}} \right).
\end{equation}
\end{theorem}

\begin{proof}
  To simplify the integration process, assume that  $A = b^2 f_1(b^2)$ and $B = f_2(b^2) - 1 = \frac{1}{4}A^2 - 1$. Substituting these into \eqref{Eq:Case_3_DE}, the denominator factors completely as a perfect square 
\begin{equation*}
\mathcal{D} = b^2 + A s \sqrt{b^2-s^2} + s^2 \left( \frac{1}{4}A^2 - 1 \right) = \left( \frac{1}{2}A s + \sqrt{b^2-s^2} \right)^2.
\end{equation*}
That is, we have 
\begin{equation}\label{eq:int_form_explicit}
\ln \phi = \int \frac{4A \sqrt{b^2-s^2} + s\left(A^2 - 4\right)}{\left( A s + 2\sqrt{b^2-s^2} \right)^2} \, ds.
\end{equation}
By decomposing the integrand fraction, we can write 
$$\ln \phi = \int\frac{A \sqrt{b^2-s^2}-2s }{\sqrt{b^2-s^2}\left( A s + 2\sqrt{b^2-s^2} \right)}\,ds+\int\frac{2b^2 A}{\sqrt{b^2-s^2}\left( A s + 2\sqrt{b^2-s^2} \right)^2}\,ds.$$
The first integral is given by
$$
\int\frac{A \sqrt{b^2-s^2}-2s }{\sqrt{b^2-s^2}\left( A s + 2\sqrt{b^2-s^2} \right)}\,ds = \ln\left(A s + 2\sqrt{b^2-s^2} \right)=\ln(b^2f_1(b^2)s + 2\sqrt{b^2-s^2}).
$$
The second one can be  evaluated as follows
$$\int\frac{2b^2 A}{\sqrt{b^2-s^2}\left( A s + 2\sqrt{b^2-s^2} \right)^2}\,ds =-\frac{2\sqrt{b^2-s^2} }{b^2 f_1(b^2) s + 2\sqrt{b^2-s^2}}.$$
Combining the above identities yields the desired expression for $\phi$.
\end{proof}

\section*{Conclusion}

In this work, we have derived explicit formulae for all possible unicorns in the settings of both standard and general $(\alpha,\beta)$-metrics, and have identified several new and easily verifiable classes of such unicorns. Furthermore, the non-degeneracy conditions \eqref{Eq:simple_form_almost regularity_1} and \eqref{Eq:simple_form_almost regularity_2} provide simple and effective criteria for testing the resulting metrics without requiring the explicit determination of the function $\phi$. To the best of our knowledge, the non-degeneracy properties of unicorn classes in the contexts of both standard and general $(\alpha,\beta)$-metrics have not been investigated previously. Moreover, the exceptional cases excluded by these conditions are themselves new and warrant further investigation, potentially opening several promising directions for future research.

\section*{Declarations}

\begin{itemize}
\item \textbf{ Declaration of Generative AI and AI-assisted technologies in the writing process}: During the preparation of this work, the author used Gemini and ChatGPT to refine the language and improve the structural flow of the manuscript. After using these tools, the author  reviewed and edited the content as needed and take full responsibility for the content of the published article.
    \item \textbf{Competing interests}: The author declares no conflict of interest.
    \item \textbf{Availability of data and material}: Not applicable.
    \item \textbf{Funding}: Not applicable.
\end{itemize}


\providecommand{\bysame}{\leavevmode\hbox to3em{\hrulefill}\thinspace}
\providecommand{\MR}{\relax\ifhmode\unskip\space\fi MR }
\providecommand{\MRhref}[2]{%
  \href{http://www.ams.org/mathscinet-getitem?mr=#1}{#2}
}
\providecommand{\href}[2]{#2}


\begin{thebibliography}{99}
\bibitem{Asa06}
G. S. Asanov, \emph{Finsleroid-Finsler spaces of Landsberg type of any dimension}, Publ. Math. Debrecen, \textbf{69} (2006), no. 1-2, 163--191.

\bibitem{Bao07}
D. Bao, \emph{On two curvature-driven problems in Finsler geometry}, Adv. Stud. Pure Math., \textbf{48} (2007), 19--71.


\bibitem{Crampin_2} 
M. Crampin, \emph{The current state of play in the Landsberg-Berwald problem of Finsler geometry}, Extracta Math., \textbf{39} (1) (2024), 57--95.


\bibitem{Elgendi-SSM} 
S. G. Elgendi, \emph{On the classification of Landsberg spherically symmetric Finsler metrics}, Int. J. Geom. Methods Mod. Phys., \textbf{18} (2021).

\bibitem{Elgendi-solutions} 
S. G. Elgendi, \emph{Solutions for the Landsberg unicorn problem in Finsler geometry}, J. Geom. Phys., \textbf{159} (2021), arXiv:1908.10910 [math.DG].

\bibitem{Elgendi-CM} 
S. G. Elgendi, \emph{New constructions of Landsberg non-Berwaldian $(\alpha, \beta)$-metrics}, Contemp. Math., (2025).

\bibitem{Elgendi-JGP} 
S. G. Elgendi, \emph{Spherically symmetric solutions to the Landsberg's unicorn problem}, J. Geom. Phys., \textbf{122} (2026).

 
\bibitem{Elgendi-MJM} 
S. G. Elgendi and Z. Muzsnay, \emph{New insights on the unicorn Landsberg problem in Finsler geometry}, Mediterr. J. Math., \textbf{22} (2025).

\bibitem{Heefer2023}
S. Heefer, C. Pfeifer, A. Reggio, and A. Fuster, \emph{A Cosmological unicorn solution to Finsler gravity}, Phys. Rev. D, \textbf{108}(6) (2023), 064051.

\bibitem{Shen2009}
Z. Shen, \emph{On a class of Landsberg metrics in Finsler geometry}, Canad. J. Math., \textbf{61} (2009), 1357--1374.



\bibitem{Tayebi2021} 
A. Tayebi, \emph{A survey on unicorns in Finsler geometry}, J. Math. Com., \textbf{2} (2) (2021), 239--250.

\bibitem{Yu-Zhu}
C. Yu and H. Zhu, \emph{On a new class of Finsler metrics}, Differ. Geom. Appl., \textbf{29} (2011), 244--554.

\bibitem{ZhouLi2018}
S. Zhou and B. Li, \emph{On Landsberg general $(\alpha, \beta)$-metrics with a conformal 1-form}, Differ. Geom. Appl., \textbf{59} (2018), 105--123.

\bibitem{Zhou-Wang-Li}
S. Zhou, J. Wang, and B. Li, \emph{On a class of almost regular Landsberg metrics}, Sci. China Math., \textbf{62} (2019), no. 5, 935--960.
 
\end{thebibliography}
\end{document}